	\documentclass{article}
	\usepackage{epsfig,latexsym}

\begin{document}

	\title{Remarks on a Ramsey theory for trees}

	\author{{\sl J\'anos Pach}\thanks{Supported by NSF Grant CCF-08-30272, and by grants from NSA, SNF, PSC-CUNY, OTKA, and BSF.} \\
	EPFL, Lausanne and R\'enyi Institute, Budapest\\
	\and {\sl J\'ozsef Solymosi}\thanks{Supported by NSERC and OTKA
	grants.} \\
	University of British Columbia, Vancouver\\
	\and {\sl G\'abor Tardos}\thanks{Supported by NSERC grant 329527
	and by OTKA grants T-046234, AT048826 and NK-62321.}\\
	Simon Fraser University and R\'enyi Institute, Budapest}

	\date{}

	\maketitle

	{\em Dedicated to Endre Szemer\'edi on the occasion of his 70th birthday.}

	\begin{abstract}
	Extending Furstenberg's ergodic theoretic proof for Szemer\'edi's theorem on arithmetic progressions, Furstenberg and Weiss (2003) proved the following qualitative result. For every $d$ and $k$, there exists an integer $N$ such that no matter how we color the vertices of a complete binary tree $T_N$ of depth $N$ with $k$ colors, we can find a monochromatic replica of $T_d$ in $T_N$ such that (1) all vertices at the same level in $T_d$ are mapped into vertices at the same level in $T_N$; (2) if a vertex $x\in V(T_d)$ is mapped into a vertex $y$ in $T_N$, then the two children of $x$ are mapped into descendants of the the two children of $y$ in $T_N$, respectively; and (3) the levels occupied by this replica form an arithmetic progression in $\{0,1,\ldots,N-1\}$. This result and its density versions imply van der Waerden's and Szemer\'edi's theorems, and laid the foundations of a new Ramsey theory for trees.

	Using simple counting arguments and a randomized coloring algorithm called {\em random split}, we prove the following related result. Let $N=N(d,k)$ denote the smallest positive integer such that no matter how we color the vertices of a complete binary tree $T_N$ of depth $N$ with $k$ colors, we can find a monochromatic replica of $T_d$ in $T_N$ which satisfies properties (1) and (2) above. Then we have $N(d,k)=\Theta(dk\log k)$. We also prove a density version of this result, which, combined with Szemer\'edi's theorem, provides a
	very short combinatorial proof of a quantitative version of the Furstenberg-Weiss theorem.
	\end{abstract}

	\section{Introduction}

	Van der Waerden's celebrated theorem \cite{vdW27} states that for any positive integers $d$ and $k$, there exists an integer $M=M(d,k)$ such that no matter how we color the elements of the set $\{1, 2,\ldots, M\}$ with $k$ colors, at least one of the color classes contains an arithmetic progression of length $d$.

	Erd\H os and Tur\'an \cite{ErT36} conjectured in 1936 and Szemer\'edi \cite{Sze75} proved in 1974 that this statement can be generalized as follows. For any positive integer $d$ and real $\delta>0$, there exists an integer $m=m(d,\delta)$ such that every subset of the set $\{1, 2,\ldots, m\}$ of size at least $\delta m$ contains an arithmetic progression of length $d$. Clearly, in van der Waerden's theorem, $M(d,k)$ can be chosen to be $m(d,1/k)$.

	A second proof of Szemer\'edi's theorem was given by Furstenberg \cite{Fu77}, using ergodic theory. Although qualitative in nature, this proof also has a quantitative version \cite{Ta06}. Furstenberg's proof represented a breakthrough, partly because of its flexibility. It led to a number of generalizations of Szemer\'edi's theorem that do not seem to follow by the original approach. These include the density Hales-Jewett theorem \cite{FuK91} and the polynomial Szemer\'edi theorem \cite{BeL96}, \cite{BeL99}.

	\smallskip

	In 2003, Furstenberg and Weiss \cite{FuW03} extended Furstenberg's proof to recurrence properties for Markov processes, which resulted in a series of new Ramsey-type theorems for trees. To formulate their results, we need to introduce some definitions.

	\smallskip

	For any positive integer $d$, let $T_d$ denote the full binary tree of depth $d-1$. We will use the terms of {root}, {leaf}, {child}, {descendant}, and {level} in their usual meaning. In the standard implementation, for any $d>0$, the vertex set $V(T_d)$ of $T_d$ consists of the strings of length smaller than $d$ over the binary alphabet $\{0,1\}$. The {\em level} of a vertex is the length of the string. The {\em root}, the only vertex at level $0$, is the empty string. The {\em leaves} are the vertices at level $d-1$. The {\em children} of a non-leaf vertex $x$ are $x0$ and $x1$. Finally, $x$ is a {\em descendant} of $y$ if $y$ is an initial segment of $x$. Any vertex is considered a descendant of itself. The empty tree will be denoted by $T_0$.

	\smallskip

	We call a function $f:V(T_d)\to V(T_n)$ a {\em regular embedding} of $T_d$ in $T_n$ if the following two conditions are satisfied.
	\begin{enumerate}
	\item If $y$ and $z$ are the two children of $x$ in $T_d$, then $f(y)$ and $f(z)$ are descendants of distinct children of $f(x)$ in $T_n$.
	\item If $x$ and $y$ are vertices at the same level of $T_d$, then $f(x)$ and $f(y)$ are also at the same level in $T_n$.
	\end{enumerate}
	For any subset $H\subseteq V(T_n)$, we say that $H$ {\em contains a
	replica} of $T_d$ if there is a regular embedding $f:V(T_d)\to H$. If, in
	addition, there exist suitable integers $a$ and $b$ such that every vertex at
	level $i$ in $T_d$ is mapped into a vertex at level $a+ib$ in $T_n$, then $H$
	is said to contain an {\em arithmetic replica} of $T_d$. See Figure~1 for an
	example of a non-arithmetic replica of $T_3$ in $T_5$.

	\bigskip
	\epsfxsize=8truecm
	\centerline{\epsffile{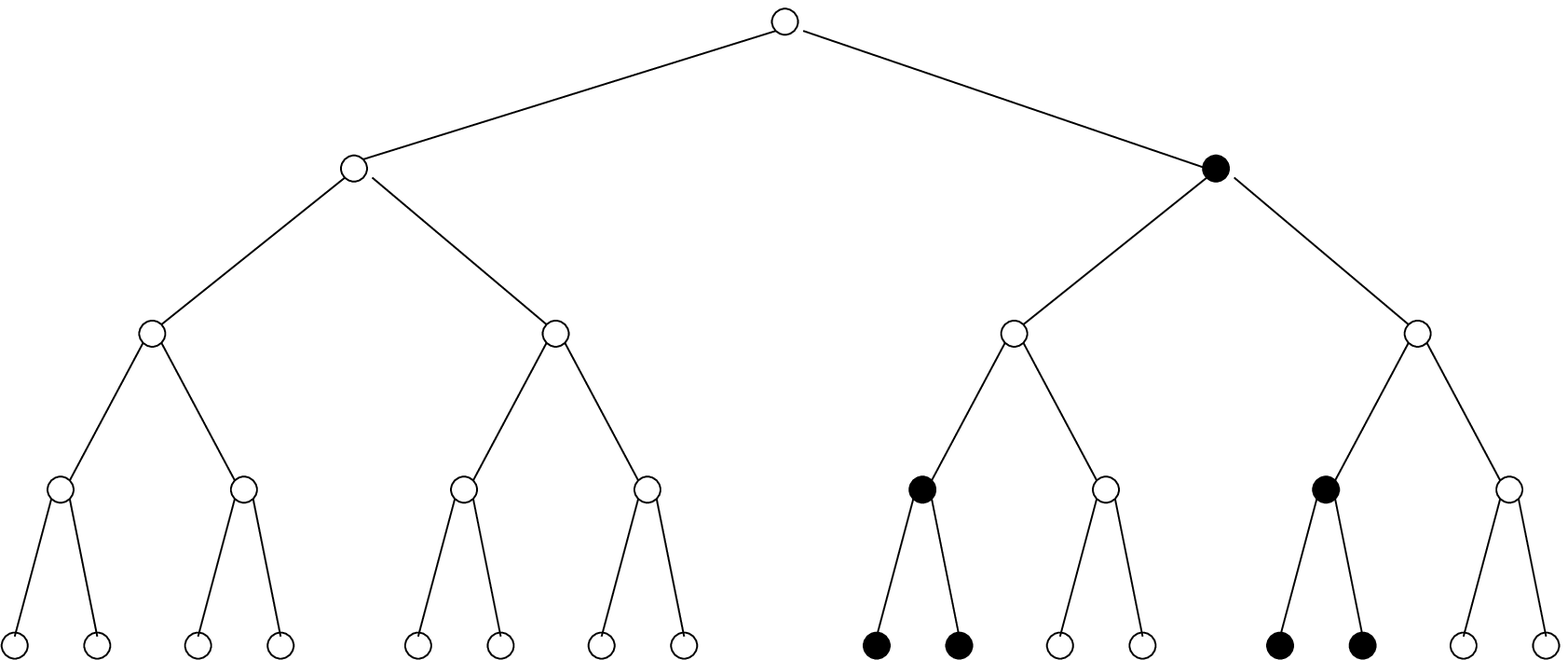}}
	\vskip0.5cm
	\centerline{{\bf Figure 1.} A non-arithmetic replica of $T_3$ in $T_5$.}
	\bigskip

	Furstenberg and Weiss \cite{FuW03} established the following theorem and various density versions of it.

	\medskip

	\noindent{\bf Theorem A.} \cite{FuW03} {\em For any positive integers $d$ and $k$, there exists $N=N(d,k)$ such that for every coloring of the vertices of $T_N$ with $k$ colors, at least one of the color classes contains an arithmetic replica of $T_d$.}

	\medskip

	Restricting this result to colorings of $T_N$, in which all vertices at the same level receive the same color, we obtain van der Waerden's theorem.

	\smallskip
	More than half of the vertices of $T_N$ are leaves, yet the set of leaves of $T_N$ contains no replica of $T_d$. Therefore, to formulate an analogue of Szemer\'edi's theorem for trees, we have to measure the ``density" of a subset $H\subset V(T_N)$ differently.

	Furstenberg and Weiss defined the {\em weight $w(x)$ of a vertex} $x\in V(T_n)$ to be $2^{-l(x)}$, where $l(x)$ denotes the level of $x$ in $T_n$. The {\em weight of a set} $H\subseteq V(T_n)$ is
	$$w(H)=\sum_{x\in H}w(x).$$
	In other words, $w(H)$ is the expected size of the intersection of $H$ with a uniformly selected random branch of $T_n$.

	\medskip

	\noindent{\bf Theorem B.} \cite{FuW03} {\em For any positive integer $d$ and real $\delta>0$, there exists $n=n(d,\delta)$ such that every subset of the vertex set of $T_n$ with weight at least $\delta n$ contains an arithmetic replica of $T_d$.}

	\medskip

	Obviously, Theorem B generalizes both Theorem A and Szemer\'edi's theorem on arithmetic progressions.

	\smallskip

	The aim of this note is to offer a simple alternative approach to Theorems A and B. Using elementary combinatorial arguments and a randomized coloring algorithm, called {\em random split}, we prove the following results.

	\bigskip

	\noindent{\bf Theorem 1.} {\em Let $d$, $n$ be positive integers, and let $H$ be a subset of the vertex set of $T_n$ satisfying
	$$2^{w(H)}>\sum_{i=0}^{d-1}{n\choose i}.$$
	Then $H$ contains a replica of $T_d$.}

	\bigskip

	\noindent{\bf Theorem 2.} {\em Let  $k, d, n\ge 2$ be integers.

	(i) Suppose that $n>5dk\log k,$ where $\log$ denotes the logarithm of base 2. Then, for any coloring of the vertices of $T_n$ with $k$ colors, one can find in $T_n$ a monochromatic replica of $T_d$.

	(ii) If $n\le dk\log k/6$, then there exists a coloring of $T_n$ with $k$ colors such that $T_n$ contains no monochromatic replica of $T_d$.}

	\medskip

	The first statement of Theorem 2 directly follows from Theorem 1. Indeed, for any $k$-coloring of $V(T_n)$, the weight of at least one of the color classes is at least $w(V(T_n))/k=n/k$. Thus, this color class contains a monochromatic replica of $T_d$, whenever we have
	$$2^{n/k}>\sum_{i=0}^{d-1}{n\choose i}.$$
	It follows by straightforward computation that this inequality holds for $n>5dk\log k.$

	\smallskip

	By its nature, the original ergodic proof of the Furstenberg-Weiss theorem is purely existential. We finish this section by showing that Theorem 1 implies a {\em quantitative} version of Theorem B with
	$$n(d,\delta)<2^{2^{(1/\delta)^{2^{2^{d+9}}}}}.$$

	\medskip
	\noindent{\bf Proof of Theorem B.} Let $H\subset V(T_n)$ be a set of weight at least $\delta n$, and let $l$ be a positive integer. We are going to prove that $H$ contains an arithmetic replica of $T_l$, provided that $n$ is sufficiently large.

	Let $d=\varepsilon n$, for some $\varepsilon>0$ to be specified later. It follows from Chernoff's bound \cite{Re70} that
	$$\sum_{i=0}^{d-1}{n\choose i}< 2^{h(\varepsilon)n},$$
	where $h(\varepsilon)=-\varepsilon\log\varepsilon-(1-\varepsilon)\log(1-\varepsilon)$ stands for the binary entropy of $\varepsilon$. Therefore, as long as $h(\varepsilon)\le\delta$, we have
	$$2^{w(H)}\ge 2^{\delta n}\ge 2^{h(\varepsilon)n}>\sum_{i=0}^{d-1}{n\choose i},$$
	and the condition in Theorem 1 is satisfied. Setting $\varepsilon=h^{-1}(\delta)$, we obtain that $T_n$ contains a replica of $T_d$ with $d=\lfloor\varepsilon n\rfloor=\lfloor h^{-1}(\delta) n\rfloor$. Let $H_d\subset H$ denote the set of elements of such a replica.

	The set of levels occupied by the elements of $H_d$ in $T_n$ is a subset of $\{0,1,\ldots, n-1\}$ with density roughly $h^{-1}(\delta)>0$. Thus, it follows from Szemer\'edi's theorem that this set contains an arithmetic progression of length $l$, provided that $n$ is sufficiently large. This implies that there is a regular embedding $f:V(T_l)\to V(H_d)$ such that the levels of $f(V(T_l))$ in $T_n$ form an arithmetic progression. In other words, $T_n$ contains an arithmetic replica of $T_l$, as desired. If we plug in the best known quantitative version of Szemer\'edi's theorem, due to Gowers \cite{Go01}, we obtain the desired bound. $\Box$

	\section{Proof of Theorem 1}

	Let us start with a couple of definitions.

	The {\it signature} of a regular embedding $f:V(T_d)\to V(T_n)$ is defined as the set of levels in $T_n$ occupied by the images of the vertices $v\in V(T_d)$. Since all vertices at the same level of $T_d$ are mapped by $f$ into vertices at the same level of $T_n$, and vertices at different levels in $V(T_d)$ are mapped into vertices at different levels, we obtain that the signature of $f$ is a $d$-element subset of $\{0,1,\ldots,n-1\}$.

	For a given subset $H\subset T_n$, We write $S(H)$ for the set
	of signatures of all regular embeddings of $T_d$ in $H$, with $d\ge0$. We have that $\emptyset\in S(H)$, which corresponds to the degenerate case when $d=0$, and $T_0$ has no vertices.

	\medskip
	\noindent{\bf Lemma 3.} {\em Let $H\subseteq V(T_n)$. We have $|S(H)|\ge2^{w(H)}$.}
	\medskip

	\noindent{\bf Proof.} The proof is by induction on $n$. If $n=0$, we have $H=\emptyset$, $w(H)=0$, and $S(H)=\{\emptyset\}$, so the statement is true.

	Suppose now that $n\ge1$ and that we have proved Lemma 3 for $n-1$. Let $r$ denote the root of $T_n$, and let $T'$ and $T''$ be the two subtrees isomorphic to $T_{n-1}$ that $T_n-r$ splits into. We apply the lemma to these subtrees and to the sets $H'=H\cap V(T')$ and $H''=H\cap V(T'')$. By the induction hypothesis, we have $|S(H')|\ge2^{2w(H')}$. Note that the weight of $H'$ inside the tree $T'$ is $2w(H')$, because the levels are shifted by one. However, this shift does not affect the {\it size} of the set of signatures. Analogously, we have $|S(H'')|\ge2^{2w(H'')}$.

	\smallskip

	We distinguish two cases. If $r\notin H$, then $w(H)=w(H')+w(H'')$ and $S(H)=S(H')\cup S(H'')$. The inequality claimed in the lemma follows:
	$$|S(H)|\ge\max(|S(H')|,|S(H'')|)\ge2^{\max(2w(H'),2w(H''))}\ge2^{w(H')+w(H'')}.$$
	On the other hand, if $r\in H$, then we have $w(H)=w(H')+w(H'')+1$, as $w(r)=1$. In view of the fact that $S(H)$ is the disjoint
	union of the sets $S(H')\cup S(H'')$ and $\{s\cup\{0\}\mid s\in S(H')\cap S(H'')\},$ we obtain that in this case $|S(H)|=|S(H')|+|S(H'')|$. Using the convexity of the function $2^w$, we can conclude that in this case
	$$|S(H)|=|S(H')|+|S(H'')|\ge2^{2w(H')}+2^{2w(H'')}\ge2\cdot2^{w(H')+w(H'')}=2^{w(H)},$$
	as desired. $\Box$

	\medskip

	\noindent{\bf Proof of Theorem 1.} By Lemma~3, the number of signatures determined by $H$ is at least $2^{w(H)}$. By the assumption of Theorem 1, this quantity is larger than the number of signatures of size smaller than $d$. Therefore, $S(H)$ has an element of size at least $d$. In other words, there exists a regular embedding of $T_d$ in $H$. $\Box$

	\section{\bf Random split and fit---Proof of Theorem 2(ii)}

	To prove the existence of a coloring which meets the requirements in Theorem~2(ii), we use a random coloring algorithm which will be called {\it random split}.

	We color the vertices of $T_n$ by the positive integers in order of
	increasing level (breadth first). While performing the coloring procedure, we maintain a list of ``forbidden colors'' for each vertex of $T_n$ not yet colored. These lists are empty at the beginning of the procedure. When we reach a vertex $x$, we assign to $x$ the smallest permitted color: the smallest positive integer $c$ that does not appear on its list of forbidden colors. If $x$ is not a leaf,
	we update the lists associated to its descendants as follows. Let $y$ and $z$ be the two children of $x$, and let $D_y$ and $D_z$ denote their sets of descendants. For each level $l$ larger than the level of $x$, we make an independent uniform random choice and either add $c$ to the list of forbidden colors of every element of $D_y$ on level $l$ or we add $c$ to the list of every vertex in $D_z$ on level $l$.

	\medskip

	\noindent{\bf Lemma 4.} {\em The random split coloring of $T_n$ admits no monochromatic regular embedding (replica) of $T_2$.}

	\medskip

	\noindent{\bf Proof.} Consider any regular embedding $f$ of $T_2$ in
	$T_n$. Denote by $x$ the image of the root of $T_2$ and let $y$ and $z$ be the images of the leaves. By definition, $y$ and $z$ are on the same level $l$ and they are descendants of distinct children of $x$. In the random split coloring, $x$ receives some color $c$ and at the same time the $c$ is added to the list of forbidden colors to all descendants of one of its children on level $l$. In particular, $c$ will be forbidden either for $y$ or for $z$. Thus, $f$ cannot be monochromatic with respect to this coloring. $\Box$

	\medskip
	\noindent {\bf Lemma 5.} {\em Restricted to any one root-leaf branch of $T_n$, the random split coloring is equivalent to the following ``random fit'' procedure: We color the vertices one by one, starting at the root. For each vertex, we consider the positive integers in increasing order until one is accepted and give the vertex the
	accepted color. When considering the integer $c$, we accept it with probability $2^{-m}$, where $m$ is the number of vertices (along this branch) that have earlier been colored with the color $c$. In particular, we do accept $c$ the first time it is considered.}

	\medskip

	\noindent{\bf Proof.} Restricting our attention to a single branch simplifies the procedure of updating the lists of forbidden colors in random split: once a vertex is colored, its color is added to the list of each uncolored vertex independently with probability $1/2$. Equivalently, in the random fit procedure, if a color $c$ appears $m$ times along the branch up to a certain point, then for every remaining vertex
	$y$, the color $c$ appears on the list of $y$ with probability $1-2^{-m}$. These events are
	independent for different pairs $(c,y)$, so deciding them can be postponed until the
	particular vertex $y$ is colored, as done in random fit. $\Box$

	\medskip

	The key to the proof of Theorem 2(ii) is the following statement.

	\medskip
	\noindent{\bf Lemma 6.} {\em Let $n\ge 8$ and $k=2\lfloor3n/\log n\rfloor$. For any branch of $T_n$ of length $n$, the probability that the random fit algorithm uses a color higher than $k$ is smaller than $2^{1-n}$.}
	\medskip

	Before proving Theorem 2(ii) in its full generality, we show that Lemma~6 implies the result for $d=2$. Indeed, using the fact that $T_n$ has $2^{n-1}$ branches, the probability that all of them will be colored by at most $k$ colors is positive. In view of Lemmas~4 and~5, this means that the coloring obtained by the random split algorithm does not admit a monochromatic replica of $T_2$, and it uses at most $k$ colors with positive probability.

	\medskip
	\noindent{\bf Proof of Lemma 6.} Fix a branch of length $n$ of $T_n$, and
	consider one by one the sequence of all choices made by the random fit
	algorithm. The maximum number of choices is $N:={n+1\choose2}$, and after each
	choice we either accept or reject a color. Let $p_j$ denote the probability
	with which we accept the color at the $j$'th choice.

	\smallskip

	Set $X_0=0$, and for any $j>0$, define the random variable $X_j$ as
	follows. Let $X_j=X_{j-1}+p_j$ if random fit rejects the corresponding color considered at the $j$'th choice, and let $X_j=X_{j-1}+p_j-1$ if it accepts. If random fit makes fewer
	than $j$ individual choices, we simply set $X_j=X_{j-1}$. Obviously, the random variables $X_j\;(j=1,2,3,...)$ define a martingale with differences bounded by $1$, and $X_j$ stabilizes for $j\ge N$.

	\smallskip

	There are exactly $n$ choices at which a color is accepted, and the corresponding $-1$ terms contribute $-n$ to $X_N$. If a color larger than $k$ was ever used, then every color up to $k+1$ must have been used at least once. For simplicity, we set $l=k/2+1$ and use the fact that each color $i\le l$ must have been considered at least $l$ times, and every time it was considered, it gave a positive contribution to $X_n$ of at least $2^{-m_i}$, where $m_i$ is the total number of vertices along this branch that were assigned color $i$. Thus, we have
	$$X_n\ge\sum_{i=1}^l\frac
	l{2^{m_i}}-n\ge\frac{l^2}{n^{1/3}}-n\ge\frac{9n^{5/3}}{\log^2n}-n,$$
	where the middle inequality comes from the fact the $\sum_i m_i\le n$ and $2^{-m}$ is a convex function.

	Azuma's inequality~\cite{AlS08} bounds the probability that $X_N=X_N-X_0>T$ by $e^{-\frac{T^2}{2N}}$. Substituting
	$T=\frac{9n^{5/3}}{\log^2n}-n$, we obtain the desired bound for the probability that a color larger than $k$ is assigned to some vertex.
	$\Box$

	\medskip

	\noindent{\bf Proof of Theorem 2(ii)}. We have already seen that for $d=2$ the statement directly follows from Lemma 6. This means that there is a $k$-coloring $\chi'$ of $T_{n'}$ with $n'=\Theta(k\log k)$, which does not admit a monochromatic regular embedding of $T_2$.

	To tackle the case $d>2$, let $n=(d-1)n'$ and split $T_n$ into subtrees isomorphic to $T_{n'}$, in the usual way: the levels 0 to $n'-1$ form one subtree, the levels $n'$ to $2n'-1$ form $2^{n'}$ subtrees, etc. Coloring each of these  subtrees separately according to $\chi'$, we obtain a coloring that admits no monochromatic regular embedding of $T_d$. $\Box$

	\section{Concluding remarks}

	Furstenberg and Weiss generalized Theorem B in two directions. First of all, instead of binary trees, one can consider ternary trees or, in general, trees in which every non-leaf vertex has $s$ children, for some $s\ge 2$.

	Obviously, our approach also applies to this case. Let $T_{n,s}$ denote a full tree of depth $n$ with this property. The only difference in our argument is that now the {\em weight} of a vertex $x\in V(T_{n,s})$ at level $l$ has to be defined as $w(x)=s^{-l}$. The weight of a subset of $V(T_{n,s})$ is the sum of the weights of its elements. We can define the regular embeddings of $T_{d,s}$ in
	$T_{n,s}$ as in the binary case, but now the $s$ children of a
	vertex $v\in V(T_{d,s})$ have to be mapped to descendants of distinct children of the image of $v$. Instead of Lemma 3, now we have

	\medskip
	\noindent{\bf Lemma 3'.} {\em Let $H\subseteq V(T_{n,s})$. Then the number of signatures of all regular embeddings of $T_{d,s}$ in $H$ satisfies $$|S(H)|\ge\sum_{\sigma\in S(H)}(s-1)^{-|\sigma|}\ge \left(\frac{s}{s-1}\right)^{w(H)}.$$}
	\medskip

	Lemma~3' readily implies the following version of Theorem~1:

	\medskip
	\noindent{\bf Theorem~1'.} {\em Let $d$, $n$, and $s$ be positive integers, and let $H$ be a subset of the vertex set of $T_{n,s}$ satisfying
	$$\left(\frac s{s-1}\right)^{w(H)}>\sum_{i=0}^{d-1}\frac{{n\choose i}}{(s-1)^i}.$$
	Then $H$ contains a replica of $T_{d,s}$.}
	\medskip

	Using Theorem~1' and Szemer\'edi's theorem, one can easily deduce the corresponding version of Theorem~B: Any subset $H$ of the vertex set of $T_{n,s}$ with weight $w(H)\ge\delta n$ contains an arithmetic replica of $T_{d,s}$ provided $n>n_0(d,s,\delta)$.

	\medskip

	As another variant of their result, Furstenberg and Weiss considered arithmetic embeddings in not necessarily full trees. Nevertheless, in what follows, we use the same definition of regular embeddings of $T_d$ in $T$ as for embeddings in $T_n$.

	\medskip
	\noindent{\bf Theorem B'.} \cite{FuW03}
	{\em Let $\alpha$ be a positive real, let $s>1$ be an integer, and $T$ a rooted tree satisfying the following three conditions:

	(a) every vertex has {\em at most} $s$ children;

	(b) every leaf is at level $n$; and

	(c) the number of leaves is at least $2^{\alpha n}$.

	\noindent Then there is an arithmetic replica of $T_d$ in $T$ provided $n>n'_0(s,d,\alpha)$.}

	\medskip

	\noindent{\bf Proof.} 
	Let us define the map $g:V(T_n)\to V(T)$ as follows. For the root $r$ of $T_n$, let $g(r)$ be the root of $T$. For any non-leaf vertex $v\in V(T_d)$, let $g$ map the two children of $v$ to the two
	distinct children of $g(v)$ which have the largest number of descendants that are leaves, unless $g(v)$ has only one child. In the latter case, we map both children of $v$ to the only child of $g(v)$.
	Let
	$$H=\{v\in V(T_d)\mid g(v)\mbox{ has 0 or at least 2 children}\}.$$
	Since $g$ preserves levels, it maps every replica of $T_d$ in $H$ into a replica of $T_d$ in $T$, and every arithmetic replica is mapped into an arithmetic replica.

	Note that $T$ has at most $s^{w(H)-1}$ leaves. Thus, by our assumption, $w(H)>\alpha n/\log s$. By Theorem~B, if $n>n'_0(s,d,\alpha)$, then $H$ contains an arithmetic replica $X$ of $T_d$. Consequently, $g(X)$ is an arithmetic replica of
	$T_d$ in $T$. $\Box$

\end{document}